\documentclass[11pt]{article}
\usepackage{amsfonts}
\usepackage{amsfonts,mathrsfs,amssymb,amsthm,mathptm}
\usepackage{amsmath,amscd}
\usepackage{mathptm,pslatex}
\usepackage{color}
\usepackage[dvips]{graphicx}
\usepackage[all]{xy}
\usepackage{graphicx}
\usepackage{fancyhdr}

\begin{document}
\newtheorem{Def}{Definition}[section]
\newtheorem{Bsp}[Def]{Example}
\newtheorem{Prop}[Def]{Proposition}
\newtheorem{Theo}[Def]{Theorem}
\newtheorem{Lem}[Def]{Lemma}
\newtheorem{Koro}[Def]{Corollary}
\theoremstyle{definition}
\newtheorem{Rem}[Def]{Remark}

\newcommand{\add}{{\rm add}}
\newcommand{\con}{{\rm con}}
\newcommand{\gd}{{\rm gl.dim}}
\newcommand{\sd}{{\rm st.dim}}
\newcommand{\sr}{{\rm sr}}
\newcommand{\dm}{{\rm dom.dim}}
\newcommand{\cdm}{{\rm codomdim}}
\newcommand{\tdim}{{\rm dim}}
\newcommand{\E}{{\rm E}}
\newcommand{\Mor}{{\rm Morph}}
\newcommand{\End}{{\rm End}}
\newcommand{\ind}{{\rm ind}}
\newcommand{\rsd}{{\rm res.dim}}
\newcommand{\rd} {{\rm rd}}
\newcommand{\ol}{\overline}
\newcommand{\overpr}{$\hfill\square$}
\newcommand{\rad}{{\rm rad}}
\newcommand{\soc}{{\rm soc}}
\renewcommand{\top}{{\rm top}}
\newcommand{\pd}{{\rm pdim}}
\newcommand{\id}{{\rm idim}}
\newcommand{\fld}{{\rm fdim}}
\newcommand{\Fac}{{\rm Fac}}
\newcommand{\Gen}{{\rm Gen}}
\newcommand{\fd} {{\rm fin.dim}}
\newcommand{\Fd} {{\rm Fin.dim}}
\newcommand{\Pf}[1]{{\mathscr P}^{<\infty}(#1)}
\newcommand{\DTr}{{\rm DTr}}
\newcommand{\cpx}[1]{#1^{\bullet}}
\newcommand{\D}[1]{{\mathscr D}(#1)}
\newcommand{\Dz}[1]{{\mathscr D}^+(#1)}
\newcommand{\Df}[1]{{\mathscr D}^-(#1)}
\newcommand{\Db}[1]{{\mathscr D}^b(#1)}
\newcommand{\A}[1]{{\mathscr A}^-(#1)}
\newcommand{\C}[1]{{\mathscr C}(#1)}
\newcommand{\Cz}[1]{{\mathscr C}^+(#1)}
\newcommand{\Cf}[1]{{\mathscr C}^-(#1)}
\newcommand{\Cb}[1]{{\mathscr C}^b(#1)}
\newcommand{\Dc}[1]{{\mathscr D}^c(#1)}
\newcommand{\K}[1]{{\mathscr K}(#1)}
\newcommand{\Kz}[1]{{\mathscr K}^+(#1)}
\newcommand{\Kf}[1]{{\mathscr  K}^-(#1)}
\newcommand{\Kb}[1]{{\mathscr K}^b(#1)}
\newcommand{\DF}[1]{{\mathscr D}_F(#1)}

\newcommand{\Kac}[1]{{\mathscr K}_{\rm ac}(#1)}
\newcommand{\Keac}[1]{{\mathscr K}_{\mbox{\rm e-ac}}(#1)}

\newcommand{\modcat}{\ensuremath{\mbox{{\rm -mod}}}}
\newcommand{\Modcat}{\ensuremath{\mbox{{\rm -Mod}}}}
\newcommand{\Spec}{{\rm Spec}}

\newcommand{\stmc}[1]{#1\mbox{{\rm -{\underline{mod}}}}}
\newcommand{\Stmc}[1]{#1\mbox{{\rm -{\underline{Mod}}}}}
\newcommand{\prj}[1]{#1\mbox{{\rm -proj}}}
\newcommand{\inj}[1]{#1\mbox{{\rm -inj}}}
\newcommand{\Prj}[1]{#1\mbox{{\rm -Proj}}}
\newcommand{\Inj}[1]{#1\mbox{{\rm -Inj}}}
\newcommand{\PI}[1]{#1\mbox{{\rm -Prinj}}}
\newcommand{\GP}[1]{#1\mbox{{\rm -GProj}}}
\newcommand{\GI}[1]{#1\mbox{{\rm -GInj}}}
\newcommand{\gp}[1]{#1\mbox{{\rm -Gproj}}}
\newcommand{\gi}[1]{#1\mbox{{\rm -Ginj}}}

\newcommand{\opp}{^{\rm op}}
\newcommand{\otimesL}{\otimes^{\rm\mathbb L}}
\newcommand{\rHom}{{\rm\mathbb R}{\rm Hom}\,}
\newcommand{\pdim}{\pd}
\newcommand{\Hom}{{\rm Hom}}
\newcommand{\Coker}{{\rm Coker}}
\newcommand{ \Ker  }{{\rm Ker}}
\newcommand{ \Cone }{{\rm Con}}
\newcommand{ \Img  }{{\rm Im}}
\newcommand{\Ext}{{\rm Ext}}
\newcommand{\StHom}{{\rm \underline{Hom}}}
\newcommand{\StEnd}{{\rm \underline{End}}}

\newcommand{\KK}{I\!\!K}

\newcommand{\gm}{{\rm _{\Gamma_M}}}
\newcommand{\gmr}{{\rm _{\Gamma_M^R}}}

\def\vez{\varepsilon}\def\bz{\bigoplus}  \def\sz {\oplus}
\def\epa{\xrightarrow} \def\inja{\hookrightarrow}

\newcommand{\lra}{\longrightarrow}
\newcommand{\llra}{\longleftarrow}
\newcommand{\lraf}[1]{\stackrel{#1}{\lra}}
\newcommand{\llaf}[1]{\stackrel{#1}{\llra}}
\newcommand{\ra}{\rightarrow}
\newcommand{\dk}{{\rm dim_{_{k}}}}

\newcommand{\holim}{{\rm Holim}}
\newcommand{\hocolim}{{\rm Hocolim}}
\newcommand{\colim}{{\rm colim\, }}
\newcommand{\limt}{{\rm lim\, }}
\newcommand{\Add}{{\rm Add }}
\newcommand{\Prod}{{\rm Prod }}
\newcommand{\Tor}{{\rm Tor}}
\newcommand{\Cogen}{{\rm Cogen}}
\newcommand{\Tria}{{\rm Tria}}
\newcommand{\Loc}{{\rm Loc}}
\newcommand{\Coloc}{{\rm Coloc}}
\newcommand{\tria}{{\rm tria}}
\newcommand{\Con}{{\rm Con}}
\newcommand{\Thick}{{\rm Thick}}
\newcommand{\thick}{{\rm thick}}
\newcommand{\Sum}{{\rm Sum}}

{\Large \bf
\begin{center}
Derived equivalences for mirror-reflective algebras
\end{center}}

\medskip\centerline{\textbf{Hongxing Chen} and \textbf{Changchang Xi}$^*$}

\renewcommand{\thefootnote}{\alph{footnote}}
\setcounter{footnote}{-1} \footnote{ $^*$ Corresponding author.
Email: xicc@cnu.edu.cn; Fax: 0086 10 68903637.}
\renewcommand{\thefootnote}{\alph{footnote}}
\setcounter{footnote}{-1} \footnote{2010 Mathematics Subject
Classification: Primary 16E35, 18E30, 16P10; Secondary 16S10, 16S50}
\renewcommand{\thefootnote}{\alph{footnote}}
\setcounter{footnote}{-1} \footnote{Keywords: Derived equivalence; Gendo-symmetric algebra; Mirror-reflective algebra; Recollement; Symmetric algebra; Tilting complex.}

\begin{abstract}
We show that the construction of mirror-reflective algebras inherits derived equivalences of gendo-symmetric algebras. More precisely, suppose $A$ and $B$ are gendo-symmetric algebras with both $Ae$ and $Bf$ faithful projective-injective left ideals generated by idempotents $e$ in $A$ and $f$ in $B$, respectively. If $A$ and $B$ are derived equivalent, then the mirror-reflective algebras of $(A,e)$ and $(B,f)$ are derived equivalent.
\end{abstract}

\section{Introduction}\label{Introduction}

Given an (associative) algebra $A$ over a commutative ring $k$, an idempotent $e$ of $A$ and an element $\lambda$ in the center of $\Lambda:=eAe$, we introduced the mirror-reflective algebra $R(A,e,\lambda)$ of $A$ at level $(e,\lambda)$ in \cite{xcf}. Roughly speaking, this algebra has the underlying $k$-module structure $A \oplus Ae\otimes_\Lambda eA$ such that $Ae\otimes_\Lambda eA$ is an ideal in $R(A,e,\lambda)$. The specialization of $R(A,e,\lambda)$ at $\lambda=e$ is called the mirror-reflective algebra of $A$ at $e$, denoted by $R(A,e)$.
In case  that $A$ is a finite-dimensional gendo-symmetric algebra over a field $k$ and $e$ is an idempotent of $A$ such that $Ae$ is a faithful and projective-injective $A$-module, the algebra $R(A,e)$ is called simply  the mirror-reflective algebra of $A$. Such a construction can be iterated and thus supplies a series of both higher Auslander algebras and recollements of derived module categories. It turns out that a new characterisation of Tachikawa's second conjecture for symmetric algebras can be formulated in terms of stratifying ideals and recollements of derived categories (see \cite{xcf}).

Our purpose of this note is to show that the construction of mirror-reflective algebras preserves derived equivalences. More precisely, we have the following.

\begin{Theo}\label{main result}
Suppose that $A$ and $B$ are finite-dimensional gendo-symmetric algebras over a field $k$ and that $_AAe$ and $_BBf$ are faithful projective-injective modules generated by idempotents $e\in A$ and $f\in B$, respectively. If $A$ and $B$ are derived equivalent, then there is an isomorphism $\sigma:Z(eAe)\to Z(fBf)$ of algebras from the center of $eAe$ to the one of $fBf$ such that, for any $\lambda\in Z(eAe)$, the mirror-reflective algebras $R(A,e,\lambda)$ and $R(B,f, (\lambda)\sigma)$ are derived equivalent.
\end{Theo}

During the course of the proof of Theorem \ref{main result}, we will give a general construction of derived equivalences of mirror-reflective algebras of arbitrary algebras at any levels in Theorem \ref{derived equivalent}. So Theorem \ref{main result} is just its consequence.

This note is sketched as follows. In Section \ref{sect2} we provide preliminaries for the proof of the main result. This includes recalling basic definitions and proving facts on derived equivalences and on mirror-reflective algebras . In Section \ref{DEDMSA} we prove Theorem \ref{derived equivalent}.

\medskip
Acknowledgements. The research work was supported partially by the National Natural Science Foundation of China (Grants 12031014 and 12122112).

\section{Preliminaries\label{sect2}}
Let $k$ denote a commutative ring with identity. All algebras in the paper are associative $k$-algebras with identity. For an algebra $A$, we denote by $A\Modcat$ the category of all left $A$-modules. Let $A\modcat$ and $A\prj$ be the full subcategories of $A\Modcat$ consisting of finitely generated $A$-modules and finitely generated projective $A$-modules, respectively.

Given an additive category $\mathcal{A}$, $\C{\mathcal{A}}$ stands for the category of all complexes $\cpx{X} = (X^i,d_X^i)$ over $\mathcal{A}$ with cochain maps as morphisms, and $\K{\mathcal{A}}$ for the homotopy category of $\C{\mathcal{A}}$. We write $\Cb{\mathcal{A}}$ and $\Kb{\mathcal{A}}$ for the full subcategories of $\C{\mathcal{A}}$ and
$\K{\mathcal{A}}$ consisting of bounded complexes over $\mathcal{A}$, respectively. When $\mathcal{A}$ is abelian, the \emph{(unbounded) derived category} of $\mathcal{A}$ is denoted by $\D{\mathcal{A}}$, which is the localization of $\K{\mathcal A}$ at all quasi-isomorphisms.

For an algebra $A$, we simply write
$\K{A}$ for $\K{A\Modcat}$ and $\D{A}$ for $\D{A\Modcat}$. Also, $A\Modcat$ is often identified with the
full subcategory of $\D{A}$ consisting of all stalk complexes concentrated in degree $0$. For an idempotent element $e$ in $A$, the category $\Kb{\add(Ae)}$ is identified with its images in $\D{A}$ under the localization functor $\K{A}\to\D{A}$.

The composition of two maps $f:X\to Y$ and $g:Y\to Z$ of sets is written as $fg$. Thus, for a map $f: X\to Y$, we write $(x)f$ for the image of $x\in X$ under $f$.

\subsection{Derived equivalences of algebras with idempotents}
In this subsection, all $k$-algebras over a commutative ring $k$ are assumed to be projective as $k$-modules. Let $A^{\rm e}:=A\otimes_kA^{\opp}$ be the enveloping algebra of an algebra $A$, and $D$ be the functor $\Hom_k(-, k)$.

We first recall the definitions of tilting complexes and derived equivalences in \cite{Rickard1, Rickard2}.

\begin{Def} Let $A$ and $B$ be algebras.

$(1)$ A complex $P\in\Kb{\prj{A}}$ is called a \emph{tilting complex} if

\quad $(i)$ $P$ is self-orthogonal, that is, $\Hom_{\Kb{\prj{A}}}(P,P[n])=0$ for any $n\neq 0$,

\quad $(ii)$  $\add(P)$ generates $\Kb{\prj{A}}$ as a triangulated category, that is, $\Kb{\prj{A}}$ is the smallest full triangulated subcateory of $\Kb{\prj{A}}$ containing $\add(P)$ and being closed under isomorphisms.

$(2)$ A complex $T\in\D{A\otimes_kB\opp}$ is called a \emph{two-sided tilting complex} if there is
a complex $T^{\vee}\in \D{B\otimes_kA\opp}$ such that
$T\otimesL_BT^\vee\simeq A$ in $\D{A^{\rm e}}$ and $T^\vee\otimesL_A T\simeq B$ in $\D{B^{\rm e}}$.
The complex $T^{\vee}$ is called the \emph{inverse} of $T$.

$(3)$  Two algebras $A$ and $B$ are said to be \emph{derived equivalent} if $\D{A}$ and $ \D{B}$ are equivalent as triangulated categories, or equivalently, $\Kb{\prj{A}}$ and $\Kb{\prj{B}}$ are equivalent as triangulated categories.
\end{Def}

Let $T$ be a two-sided tilting complex in $\D{A\otimes_kB\opp}$ with the inverse  $T^{\vee}$.  By \cite[Section 3]{Rickard2}, we have
$T^{\vee}\simeq\rHom_A(T, A)\simeq\rHom_{B\opp}(T,B)$ in $\D{B\otimes_kA\opp}$. Moreover, the functor $T^{\vee}\otimesL_A-:\D{A}\to\D{B}$ is a triangle equivalence with the quasi-inverse $T\otimesL_B-: \D{B}\to\D{A}$.
This implies that $_AT$ and $T_B$ are isomorphic to tilting complexes in $\D{A}$ and $\D{B\opp}$, respectively. By \cite[Lemma 4.3]{Rickard2},
$T^\vee\otimes_kT\in\D{A^{\rm e}\otimes_k(B^{\rm e})\opp}$ is a two-sided tilting complex.

The following theorem is well known (see \cite{Happel, keller, Rickard1, Rickard2}).

\begin{Theo} \label{Derived equivalence}
Let $A$ and $B$ be $k$-algebras. The following are equivalent.

$(1)$ $A$ and $B$ are derived equivalent.

$(2)$ There is a tilting complex $P\in\Kb{\prj{A}}$ such that
$B\simeq\End_{\D A}(P)$ as algebras.

$(3)$ There is a two-sided tilting complex $T\in\D{A\otimes_kB\opp}$.
\end{Theo}

Comparing with recollement-tilting complexes related to idempotents in \cite[Definition 3.6]{Miyachi}, we introduce the definition of derived equivalences of algebras with idempotents.

\begin{Def}\label{IDR}
Let $A$ and $B$ be algebras with idempotent elements $e=e^2\in A$ and $f=f^2\in B$. The pairs $(A, e)$ and $(B,f)$ of algebras with idempotents
are said to be \emph{derived equivalent} provided that there is a triangle equivalence $\D{A}\to \D{B}$ which restricts to an equivalence $\Kb{\add(Ae)}\to\Kb{\add(Bf)}$.
\end{Def}

Clearly, $A$ and $B$ are derived equivalent if and only if so are the pairs $(A,0)$ and $(B,0)$ if and only if so are the pairs $(A,1_A)$ and $(B,1_B)$.
The following result is essentially implied in \cite{Miyachi} and provides several equivalent characterizations of derived equivalences of algebras with idempotents. For the convenience of the reader, we provide a proof.

\begin{Lem}{\rm (\cite{Miyachi})}\label{IMDE}
Let $A$ and $B$ be algebras with $e^2=e\in A$ and $f^2=f\in B$. The following are equivalent.

$(1)$ The pairs $(A, e)$ and $(B,f)$ are derived equivalent.

$(2)$ There is a tilting complex $P\in\Kb{\prj{A}}$ such that $P=P_1\oplus P_2$ in $\Kb{\prj{A}}$
 satisfying

\quad $(a)$ $B\simeq\End_{\D A}(P)$ as algebras.

\quad $(b)$ $P_1$ generates $\Kb{\add(Ae)}$ as a triangulated category.

\quad $(c)$ Under the isomorphism of $(a)$, $f\in B$ corresponds to the composite of the canonical projection $P\to P_1$ with the
canonical inclusion $P_1\to P$.

$(3)$ There is a two-sided tilting complex $T\in\D{A\otimes_kB\opp}$ with the inverse  $T^{\vee}\in\D{B\otimes_kA\opp}$ such that $eTf\in \D{eAe\otimes_k (fBf)\opp}$ is a two-sided tilting complex with the inverse  $fT^{\vee}e\in\D{fBf\otimes_k (eAe)\opp}$ and that all $3$ squares in the following diagram are commutative (up to natural isomorphism):
$$\xymatrix{\mathscr{D}(A)\ar[d]_-{F_1}\ar[rr]^-{e\cdot}&&\mathscr{D}(eAe)\ar[d]^-{F_2}\ar@/^1.2pc/[ll]_-{j_{e*}}\ar@/_1.2pc/[ll]_-{j_{e!}}\\
\mathscr{D}(B)\ar[rr]^-{f\cdot}&&\mathscr{D}(fBf)\ar@/^1.2pc/[ll]_-{j_{f*}}\ar@/_1.2pc/[ll]_-{j_{f!}}}$$
where $F_1:=T^\vee\otimesL_A-, \; F_2:=fT^\vee e\otimesL_{eAe}-,\; j_{e!}:=Ae\otimesL_{eAe}-,$
$j_{e*}:=\rHom_{eAe}(eA,-),\;\; j_{f!}:=Be\otimesL_{fBf}-,\; j_{f*}:=\rHom_{fBf}(fB, -),$
and the functors $e\cdot$ and $f\cdot$ denote the left multiplications by $e$ and $f$, respectively.

$(4)$ There is a two-sided tilting complex $T\in\D{A\otimes_kB\opp}$ with the inverse  $T^{\vee}\in\D{B\otimes_kA\opp}$ such that
$$T^\vee\otimesL_A(Ae\otimesL_{eAe}eA)\otimesL_AT\simeq Bf\otimesL_{fBf}fB\in\D{B^{\rm e}}.$$
\end{Lem}

{\it Proof.} $(1)\Rightarrow (2)$. Assume  $(1)$ holds. Then
there is a triangle equivalence $F_1:\D{A}\to\D{B}$ which restricts to an equivalence $\Kb{\add(Ae)}\to\Kb{\add(Bf)}$.
Let $G_1:\D{B}\to D(A)$ be the inverse of $F_1$. Define $P:=G_1(B)$, $P_1:=G(Bf)$ and $P_2:=G(B(1-f))$. Then $P=P_1\oplus P_2$ and $P_1\in\Kb{\add(Ae)}$.
Since $Bf$ generates $\Kb{\add(Bf)}$ as a triangulated category, all conditions $(a),(b)$ and $(c)$ hold.

$(2)\Rightarrow (3)$. Let $\Lambda:=eAe$. Recall that the adjoint pair $(Ae\otimes_\Lambda-, e\cdot)$ between $\Lambda\Modcat$ and $A\Modcat$
induces a triangle equivalence $\Kb{\add(Ae)}\lraf{\simeq}\Kb{\prj{\Lambda}}$. Since $P_1$ is a direct summand of $P$ and generates $\Kb{\add(Ae)}$ as a triangulated category, the complex $eP_1\in\Kb{\prj{\Lambda}}$ is a tilting complex. Let $T$ be a  two-sided tilting complex in $\D{A\otimes_kB\opp}$ which is induced by ${_A}P$.
Then the argument in the proof of \cite[Theorem 3.5]{Miyachi} shows that $(2)$ implies $(3)$.

$(3)\Rightarrow (1)$. Let $\Gamma:=fBf$. Note that the image of the restriction of $j_{e!}$ to $\Kb{\prj{\Lambda}}$ coincides with the image of $\Kb{\add(Ae)}$ in $\D{A}$.
Similarly, the image in $\D{B}$ of the restriction of $j_{f!}$ to $\Kb{\prj{\Gamma}}$ coincides with the image of $\Kb{\add(Bf)}$ in $\D{B}$.
Thus the equivalence $F_1$ in $(3)$ restricts to an equivalence from $\Kb{\add(Ae)}$ to $\Kb{\add(Bf)}$. Thus $(1)$ holds.

$(3)\Rightarrow (4)$. By \cite[Corollaries 3.7 and 3.8]{Miyachi},  there are isomorphisms in $\D{A\otimes_kB\opp}$:
$$T\otimesL_BBf\otimesL_\Gamma fB\simeq Tf\otimesL_\Gamma fB\simeq Ae\otimesL_\Lambda eT \simeq Ae\otimesL_\Lambda eA\otimesL_AT.$$
Applying $T^\vee\otimesL_A-:\D{A\otimes_kB\opp}\to \D{B^{\rm e}}$ to these isomorphisms yields
$$Bf\otimesL_\Gamma fB\simeq T^\vee\otimesL_AT\otimesL_BBf\otimesL_\Gamma fB\simeq T^\vee\otimesL_AAe\otimesL_\Lambda eA\otimesL_AT.$$

$(4)\Rightarrow (2)$.  Since ${_A}T_B$ is a two-sided tilting complex, it follows from $(4)$ that there are isomorphisms of complexes
$$\begin{array}{ll} (i) & Ae\otimesL_\Lambda eA\simeq T\otimesL_BBf\otimesL_\Gamma fB\otimesL_BT^\vee \in\D{A^{\rm e}},\\
(ii) &  Ae\otimesL_\Lambda eTf\simeq T\otimesL_BBf\otimesL_\Gamma fB\otimesL_BT^\vee\otimesL_AT\otimesL_BBf \simeq  Tf\in\D{A\otimes \Gamma\opp},\\
(iii) & fT^\vee e\otimesL_\Lambda eTf  \simeq fT^\vee\otimesL_AAe\otimesL_\Lambda eA\otimesL_ATf\simeq fBf\otimesL_\Gamma fBf\simeq\Gamma\in\D{\Gamma^{\rm e}},\\
(iv) & eTf\otimesL_\Gamma fT^\vee e\simeq eT\otimesL_BBf\otimesL_\Gamma fB\otimesL_BT^\vee e\simeq eAe\otimesL_\Lambda eAe\simeq\Lambda \in \D{\Lambda^{\rm e}}.\end{array}$$
Due to $(iii)$ and $(iv)$, ${_\Lambda}(eTf)_\Gamma$ is a two-sided tilting complex with the inverse  $fT^{\vee}e$. In particular,  ${_\Lambda}eTf$ is isomorphic to a tilting complex.
Since  $j_{e!} $ induces a triangle equivalence  $\Kb{\prj{\Lambda}}\lraf{\simeq}\Kb{\add(Ae)}$, the isomorphisms in $(ii)$ imply that $Tf$ generates $\Kb{\add(Ae)}$ as
a triangulated category.  Clearly,  ${_A}T$ is isomorphic to a tilting complex and has a direct summand $Tf$. Moreover, $\End_{\D{A}}(T)\simeq B$ as algebras and $\Hom_{\D A}(T, Tf)\simeq
Bf$ as $B$-modules. Thus $(2)$ holds. $\square$

\begin{Koro}\label{EVA}
Assume that the pairs $(A,e)$ and $(B,f)$  are derived equivalent. Then

$(1)$ $(A\opp,  e\opp)$ and $(B\opp, f\opp)$  are derived equivalent.

$(2)$ $(A^{\rm e}, e\otimes e\opp)$ and $(B^{\rm e}, f\otimes f\opp)$ are derived equivalent.
\end{Koro}

{\it Proof.}  Let $(-)^*:=\Hom_A(-, A)$ and $P$ be the tilting complex in  Lemma \ref{IMDE}(2). Then  $P^*\in\Kb{\prj{A\opp}}$ and $P^*=P_1^*\oplus P_2^*$. By \cite[Proposition 9.1]{Rickard1}, $P^*$ is a tilting complex over $A\opp$.

$(1)$ Since $(-)^*:\Kb{\prj{A}}\to \Kb{\prj{A\opp}}$  is a triangle equivalence sending $Ae$ to  $eA$, it follows from  Lemma \ref{IMDE}$(c)$ that
$P_1^*$ generates $\Kb{\add(eA)}$ as a triangulated category. By Lemma \ref{IMDE}$(a)$ and $(c)$, there is an algebra isomorphism $B\opp\simeq \End_{\D{A\opp}}(P^*)$
under which $f\opp$ is the composition of the projection $P^*\to P_1^*$ with the inclusion $ P_1^*\to P^*$.  Thus $(P^*, e\opp)$ satisfies Lemma \ref{IMDE}(2).
This shows $(1)$.

$(2)$ Let $Q:=P\otimes_kP^*\in\Kb{\prj{A^{\rm e}}}$. We will show that $Q$ satisfies Lemma \ref{IMDE}(2) for the pair $(A^{\rm e}, e\otimes e\opp)$ and $(B^{\rm e}, f\otimes f\opp)$.

In fact, by \cite[Theorem 2.1]{Rickard2}, $Q$ is a tilting complex over $A^{\rm e}$ and  $\End_{\D{A^{\rm e}}}(Q)\simeq B^{\rm e}$. Clearly, $P_1\otimes_kP_1^*$ is a direct summand of $P\otimes_kP^*$ and  there are canonical isomorphisms
$$\Hom_{\D{A^{\rm e}}}(Q, P_1\otimes_kP_1^*)\simeq \Hom_{\D{A}}(P, P_1)\otimes_k\Hom_{\D{A\opp}}(P^*,P_1^*)\simeq  Bf\otimes_k fB=B^{\rm e}(f\otimes f\opp).$$
Thus $Q$ satisfies Lemma \ref{IMDE}$(a)$-$(b)$. To show  Lemma \ref{IMDE}$(c)$ for $Q$, we need the following general result:

If $L:\mathcal{C}\to\mathcal{D}$ is a triangle functor between triangulated categories $\mathcal{C}$ and $\mathcal{D}$, then $L(\tria_\mathcal{C}(\add(X)))\subseteq\tria_\mathcal{D}(\add(L(X)))$ for any $X\in\mathcal{C}$, where $\tria_\mathcal{C}(\add(X))$ denotes the smallest full triangulated subcategory of $\mathcal{C}$ containing $\add(X)$.

Since $eA\in \Kb{\add(eA)}=\tria_{\K{A\opp}}(\add(P_1^*))$,  we apply the functor $Ae\otimes_k-: \Kb{\add(eA)}\to \Kb{\add(Ae\otimes_keA)}$ to the $k$-module $eA$ and obtain $Ae\otimes_keA\in\tria_{\K{A^{\rm e}}}(\add(Ae\otimes_kP_1^*))$. Similarly, we have $Ae\otimes P_1^*\in\tria_{\K{A^{\rm e}}}(\add (P_1\otimes_kP_1^*))$ by the functor $-\otimes_kP_1^*: \Kb{\add(Ae)}\to \Kb{\add(Ae\otimes_keA)}$.
Thus  $Ae\otimes_keA\in\tria(\add(P_1\otimes_kP_1^*))$.  By the equivalences of  Lemma \ref{IMDE}(1)-(2), the pairs $(A^{\rm e}, e\otimes e\opp)$ and $(B^{\rm e}, f\otimes f\opp)$ are derived equivalent. $\square$

\medskip
A finite-dimensional algebra $A$ over a field $k$ is called a \emph{gendo-symmetric} algebra
if $A=\End_\Lambda(\Lambda\oplus M)$ with $\Lambda$ a symmetric algebra and $M$ a finite-dimensional $\Lambda$-module.
By \cite[Theorem 3.2]{FK11}, $A$ is gendo-symmetric if and only if the dominant dimension of $A$ is at least $2$ and $D(Ae)\simeq eA$ as $eAe$-$A$-bimodules, where $e\in A$ is an idempotent element such that $Ae$  is a faithful projective-injective $A$-module.

\begin{Prop}{\rm \cite[Proposition 3.9]{FM19}}\label{DGS}
Suppose that  $A$ and $B$ are gendo-symmetric algebras with $Ae$ and $Bf$ faithful projective-injective modules over $A$ and $B$, respectively. If $A$ and $B$ are derived equivalent, then the pairs $(A,e)$ and $(B,f)$ are derived equivalent of algebras with idempotents.
\end{Prop}

\subsection{Mirror-reflective algebras}

In this section, we recall the construction of mirror-reflective algebras in \cite{xcf}. Assume that $A$ is a $k$-algebra over a commutative ring $k$, $e=e^2\in A$, $\Lambda:=eAe\;\;\mbox{and}\;\; \lambda$ lies in the center $Z(\Lambda)$ of $\Lambda$. Recall that the mirror-reflective algebra $R(A,e,\lambda)$ of $A$ at level $(e,\lambda)$, defined in \cite{xcf},  has the underlying $k$-module $A \oplus Ae\otimes_\Lambda eA$ as its abelian group. Its multiplication $\ast$ is given explicitly by
$$ (a + be\otimes ec)\ast (a' + b'e\otimes ec'):= aa' + (ab'e\otimes ec'+ be\otimes eca'+ becb'e\otimes\lambda ec')$$
for $a,b,c,a',b',c'\in A$.
This can be reformulated as follows: Let $\omega_\lambda$ be the composite of the natural maps:
{\small
$$(Ae\otimes_\Lambda eA)\otimes_A(Ae\otimes_\Lambda eA)\lraf{\simeq}Ae\otimes_\Lambda (eA\otimes_AAe)\otimes_{\Lambda} eA \lraf{\simeq} Ae\otimes_\Lambda \Lambda\otimes_\Lambda eA\lraf{{\rm Id}\otimes(\cdot\lambda)\otimes {\rm Id}}
Ae\otimes_\Lambda \Lambda\otimes_\Lambda eA\ra Ae\otimes_\Lambda eA,
$$}
where $(\cdot\lambda):\Lambda\to \Lambda$ is the multiplication map by $\lambda$.
Then $$\big((be\otimes ec)\otimes(b'e\otimes ec')\big)\omega_\lambda=(be\otimes ec)\ast (b'e\otimes ec').$$
Clearly, $R(A,e,0)$ is exactly the trivial extension of $A$ by $Ae\otimes_\Lambda eA$. To understand $R(A,e,\lambda)$, we will employ idealized extensions of algebras.

\begin{Def}\label{def2.7}
Let $X$ be an $A$-$A$-bimodule. An \emph{idealized extension} of $A$ by $X$ is defined to be an algebra $R$ such that $A$ is a subalgebra (with the same identity) of $R$, $X$ is an ideal of $R$, and $R=A\oplus X$ as $A$-$A$-bimodules.  Two idealized extensions $R_1$ and $R_2$ of $A$ by $X$ are said to be \emph{isomorphic} if there exists an algebra isomorphism $\phi: R_1\to R_2$ such that
the restriction of $\phi$ to $A$ is the identity map of $A$ and the one of $\phi$ to $X$ is an bijection from $X$ to $X$.
\end{Def}

Clearly, an algebra $R$ is an idealized extension of $A$ by $X$ if and only if $R$ contains $A$ as a subalgebra and
there is an algebra homomorphism $\pi:R\to A$ with $X=\Ker(\pi)$ such that the composite of the inclusion $A\to R$ with $\pi$ is the identity map of $A$. Hence a mirror-reflective algebra $R(A,e, \lambda)$ is an idealized extension of $A$ by $Ae\otimes_{\Lambda}eA$.

Let $$F:=Ae\otimes_\Lambda-\otimes_\Lambda eA:\;\;\Lambda^{\rm{e}}\Modcat\lra A^{\rm{e}}\Modcat, \;\;M\mapsto Ae\otimes_\Lambda M\otimes_\Lambda eA,$$
$$G:=e(-)e:\;\; A^{\rm{e}}\Modcat\lra\Lambda^{\rm{e}}\Modcat, \;\;M\mapsto eMe$$
for $M\in A^{\rm{e}}\Modcat$. Since $e\otimes e\opp$ is an idempotent element of $A^{\rm e}$ and there are natural isomorphisms
$$ F\simeq A^{\rm{e}}(e\otimes e\opp)\otimes_{\Lambda^{\rm{e}}}-\quad\mbox{and}\quad G\simeq \Hom_{A^{\rm{e}}}(A^{\rm{e}}(e\otimes e\opp), -),$$
$(F,G)$ is an adjoint pair and $F$ is fully faithful. This implies the following.

\begin{Lem}\label{iso}
The functor $F$ induces an algebra isomorphism
$$\rho:\;\; Z(\Lambda)\lra \End_{A^{\rm e}}(Ae\otimes_\Lambda eA),\;\; \lambda\mapsto \rho_\lambda:=[ae\otimes eb\mapsto ae\lambda\otimes eb]$$
for $\lambda\in Z(\Lambda)$ and $a, b\in A$. Moreover, $\omega_\lambda=\omega_e\rho_\lambda$.
\end{Lem}

The following result parameterizes the idealized extensions of $A$ by $Ae\otimes_\Lambda eA$.
\begin{Prop}
Let $Z(\Lambda)^{\times}$ be the group of units of $Z(\Lambda)$, that is, $Z(\Lambda)^{\times}$  is the group of all invertible elements in $Z(\Lambda)$. Then there exists a bijection from the quotient of the multiplicative semigroup $Z(\Lambda)$ modulo  $Z(\Lambda)^{\times}$ to the set $\mathscr{S}(A,e)$ of the isomorphism classes of idealized extensions of $A$ by $Ae\otimes_\Lambda eA$:
$$Z(\Lambda)/Z(\Lambda)^{\times}\lraf{\simeq} \mathscr{S}(A,e),\quad \lambda\, Z(\Lambda)^{\times}\mapsto R(A,e,\lambda) \; \mbox{ for } \lambda\in Z(\Lambda).$$
\end{Prop}

{\it Proof.} Let $Z_0(\Lambda):=Z(\Lambda)/Z(\Lambda)^{\times}=\{\lambda\, Z(\Lambda)^{\times}\mid \lambda\in Z(\Lambda)\}$ and $[\lambda]:=\lambda\, Z(\Lambda)^{\times}\in Z_0(\Lambda)$ for $\lambda\in Z(\Lambda)$.
By \cite[Lemma 3.2(2)]{xcf}, if $\mu\in Z(\Lambda)^{\times}$, then $R(A,e,\lambda)\simeq R(A,e,\lambda\mu)$ as algebras. This means that
the map $$\varphi: Z_0(\Lambda)\lra \mathscr{S}(A,e), [\lambda]\mapsto R(A,e,\lambda)$$
is well defined.  Let $R$ be an idealized extension of $A$ by $X:=Ae\otimes_\Lambda eA$. Then the multiplication of $R$ induces a homomorphism $\phi:X\otimes_AX\to X$ of $A^{\rm e}$-modules. Recall that $\omega_e: X\otimes_AX\to X$ is an isomorphism of $A^{\rm e}$-modules. Let $\phi':=\omega_e^{-1}\phi$. Then $\phi'\in\End_{A^{\rm e}}(X)$ and
$\phi=\omega_e\phi'$. By Lemma \ref{iso}, $\phi'=\rho_z$ for some $z\in Z(\Lambda)$ and $\phi=\omega_z$. Thus $R=R(A,e,z)$ and $\varphi$ is surjective.

Now, we show that $\varphi$ is injective. Suppose $\lambda_i\in Z(\Lambda)$ for $i=1,2$ and $R(A,e,\lambda_1)\simeq R(A,e,\lambda_2)$ as algebras. Set $R_i:=R(A,e,\lambda_i)$. By Definition \ref{def2.7}, there is an algebra isomorphism $f: R_1\to R_2$ such that $f|_A=\rm{Id}_A$  and $\alpha:=f|_X:X\to X$
is an isomorphism of ideals. This implies that $\alpha$ is a homomorphism of $A^{\rm e}$-modules and
$(\alpha\otimes_A\alpha)\omega_{\lambda_2}=\omega_{\lambda_1}\alpha: X\otimes_AX\to X$. Since $\omega_{\lambda_i}=\omega_e\rho_{\lambda_i}$ by Lemma \ref{iso}, there holds $(\alpha\otimes_A\alpha)\omega_e\rho_{\lambda_2}=\omega_e\rho_{\lambda_1}\alpha$.
Let $\sigma:=\omega_e^{-1}(\alpha\otimes_A\alpha)\omega_e\in\End_{A^{\rm e}}(X)$. Then $\sigma$ is an isomorphism of $A^{\rm e}$-modules and $\rho_{\lambda_1}\alpha=\sigma\rho_{\lambda_2}$.
Again by Lemma \ref{iso}, $\alpha=\rho_c$ and $\sigma=\rho_d$ for some $c, d\in Z(\Lambda)^{\times}$.
It follows that $\lambda_1c=d\lambda_2$, and therefore $[\lambda_1]=[\lambda_2]$. $\square$

\section{Derived equivalences of mirror-reflective algebras}\label{DEDMSA}

In this section, $k$ denotes a commutative ring, all algebras are $k$-algebras which are projective as $k$-modules.

Assume that the pairs $(A, e)$ and $(B,f)$ of algebras with idempotents are derived equivalent. By Lemma \ref{IMDE}, there is a two-sided tilting complex $T\in\D{A\otimes_kB\opp}$ with the quasi-inverse $T^\vee$ such that $T\otimes_kT^\vee\in \D{A^{\rm e}\otimes_k(B^{\rm e})\opp}$ is a two-sided tilting complex with the inverse  $T^\vee\otimes_kT$, and there is a derived equivalence:
$$\Phi:=T^\vee\otimesL_A-\otimesL_AT\simeq (T^\vee\otimes_kT)\otimesL_{A^{\rm e}}-:  \D{A^{\rm e}}\lra  \D{B^{\rm e}}$$
which sends $A$ to $B$ up to isomorphism (see \cite{Rickard2}). Let $\varepsilon_A:T\otimesL_BT^{\vee}\to A$ and $\varepsilon_B:T^{\vee}\otimesL_AT\to B$
be the associated isomorphisms in $\D{A^{\rm e}}$ and $\D{B^{\rm e}}$, respectively. Now, we introduce the notation
$$\Lambda=eAe,\;\;\Gamma=fBf,\;\;G_e=e(-)e,\;\;G_f=f(-)f,$$
$$F_e=Ae\otimes_\Lambda-\otimes_\Lambda eA: \;\Lambda^{\rm e}\Modcat\lra A^{\rm e}\Modcat,\quad F_f=Bf\otimes_\Gamma-\otimes_\Gamma fB:\; \Gamma^{\rm e}\Modcat\lra B^{\rm e}\Modcat,$$
$$\mathbb{L}F_e=Ae\otimesL_\Lambda-\otimesL_\Lambda eA: \;\D{\Lambda^{\rm e}}\lra \D{A^{\rm e}},\quad
\mathbb{L}F_f=Bf\otimesL_\Gamma-\otimesL_\Gamma fB: \;\D{\Gamma^{\rm e}}\lra \D{B^{\rm e}},
$$
$$\Phi'=fT^\vee e\otimesL_\Lambda-\otimesL_\Lambda eTf: \; \D{\Lambda^{\rm e}}\lra \D{\Gamma^{\rm e}},$$
$$
\Delta_0=Ae\otimes_\Lambda eA,\;\;\Delta=Ae\otimesL_\Lambda eA,\quad\Theta_0=Bf\otimes_\Gamma fB,\quad \Theta=Bf\otimesL_\Gamma fB,
$$
together with the identifications (up to isomorphism):
$$
\Delta_0=H^0(\Delta),\;\; \Theta_0=H^0(\Theta),\;\;\Delta=\mathbb{L}F_e(\Lambda),\;\;
\Theta=\mathbb{L}F_f(\Gamma).
$$
By Lemma \ref{IMDE} and Corollary \ref{EVA}, up to natural isomorphism, two squares in the diagram are commutative:
$$(\sharp)\quad\xymatrix{\mathscr{D}(A^{\rm e})\ar[d]_-{\Phi}\ar[rr]^-{G_e}&&\mathscr{D}(\Lambda^{\rm e})\ar[d]^-{\Phi'}\ar@/_1.2pc/[ll]_-{\mathbb{L}F_e}\\
\mathscr{D}(B^{\rm e})\ar[rr]^-{G_f}&&\mathscr{D}(\Gamma^{\rm e})\ar@/_1.2pc/[ll]_-{\mathbb{L}F_f}}$$
where $\Phi'$ is the derived equivalence associated with the two-sided tilting complex $eTf\in\D{\Lambda\otimes_k\Gamma}$.
Note that $\Phi, \Phi', \mathbb{L}F_e$ and $\mathbb{L}F_f$ commute with derived tensor products. Namely, for $U,V\in\D{A^{\rm e}}$,  there are isomorphisms
$$\Phi(U\otimesL_AV)\simeq T^\vee\otimesL_AU\otimesL_AA\otimesL_AV\otimesL_AT\simeq T^\vee\otimesL_AU\otimesL_AT\otimesL_BT^{\vee}\otimesL_AV\otimesL_AT= \Phi(U)\otimesL_B\Phi(V)$$
where the second isomorphism follows from $A\simeq T\otimesL_BT^{\vee}$ in $\D{A^{\rm e}}$. This provides a natural isomorphism
$$\phi_{-,-}:\;\;\Phi(-)\otimesL_B\Phi(-)\lraf{\simeq}\Phi(-\otimesL_A-):\quad \D{A^{\rm e}}\times \D{A^{\rm e}}\lra \D{B^{\rm e}}.$$
Since $\Phi'(\Lambda)\simeq \Gamma$, there is an  algebra isomorphism
$$\sigma: Z(\Lambda)\lra Z(\Gamma)$$
defined by the series of isomorphisms
$Z(\Lambda)\simeq \End_{\Lambda^{\rm e}}(\Lambda)\lraf{\simeq}  \End_{\Gamma^{\rm e}}(\Phi'(\Lambda))\lraf{\simeq}\End_{\Gamma^{\rm e}}(\Gamma)\simeq Z(\Gamma).$

Our main result on derived equivalences of mirror-reflective algebras is the following.

\begin{Theo}\label{derived equivalent}
Suppose that there is a derived equivalence between $(A, e)$ and $(B,f)$ of algebras with idempotents, which gives rise to a two-sided tilting complex ${_A}T_{B}$. If the derived equivalence $\Phi:\D{A^{\rm e}}\to \D{B^{\rm e}}$ associated with $T$ between the enveloping algebras $A^{\rm e}$ and $B^{\rm e}$ satisfies
$\Phi(Ae\otimes_\Lambda eA)\simeq Bf\otimes_\Gamma fB$ in $\D{B^{\rm e}}$, then there is an algebra isomorphism $\sigma: Z(\Lambda)\to Z(\Gamma)$ such that, for each $\lambda\in Z(\Lambda)$, the pairs $\big(R(A,e,\lambda), e\big)$ and $\big(R(B,f,(\lambda)\sigma), f\big)$ of algebras with idempotents are derived equivalent. In particular, $R(A,e)$ and $R(B,f)$ are derived equivalent.
\end{Theo}

Before starting with the proof of Theorem \ref{derived equivalent}, we first fix notation on derived categories.

Let $\mathcal{A}$ be an abelian category. For each $X:=(X^i, d_X^i)_{i\in\mathbb{Z}}\in\C{\mathcal{A}}$ and $n\in\mathbb{Z}$, there are two truncated complexes
$$
\tau^{\leq n}X:\;\; \cdots\lra X^{n-3}\epa{d_X^{n-3}} X^{n-2}\epa{d_X^{n-2}}X^{n-1}\epa{d_X^{n-1}}\Ker(d_X^n)\lra 0,
$$
$$
\tau^{\geq n}X:\;\;0\lra \Coker(d_X^{n-1})\epa{\overline{d_X^n}} X^{n+1}\lraf{d_X^{n+1}} X^{n+2}\lraf{d_X^{n+2}} X^{n+3}\lra \cdots,
$$
where $\overline{d_X^n}$ is induced from $d_X^n$. Moreover, there are canonical chain maps in $\C{\mathcal{A}}$:
$$\lambda_{X}^n:\tau^{\leq n}X\hookrightarrow X \;\;\mbox{and}\;\;\pi_{X}^n: X\twoheadrightarrow\tau^{\geq n}X,$$
and a distinguished triangle in $\D{\mathcal{A}}$:
$$\tau^{\leq n}X\lraf{\lambda_{X}^n}X\lraf{\pi_{X}^{n+1}}\tau^{\geq {n+1}}X\lra \tau^{\leq n}X[1].$$
Note that $H^n(X)=\tau^{\geq n}\tau^{\leq n}X: \D{\mathcal{A}}\to \mathcal{A}$. Let
$\mathscr{D}^{\leq 0}(\mathcal{A}):=\{X\in\D{\mathcal{A}}\mid H^i(X)=0, i>0\}.$
For each $X\in \mathscr{D}^{\leq 0}(\mathcal{A})$, it is clear that $\lambda_{X}^0$ is an isomorphism in $\D{\mathcal{A}}$.
In this case, we denote by $\xi_X: X\to H^0(X)$ the composition of the inverse $X\to\tau^{\leq 0}X$ of $\lambda_{X}^0$ with $\pi_{\tau^{\leq 0}X}^0:\tau^{\leq 0}X\to H^0(X)$.
Clearly, if $X^i=0$ for all $i\geq 1$, then $X=\tau^{\leq 0}X$ and $\xi_X=\pi^0_X$. Now, there is a natural transformation $$\xi:\;\;{\rm Id}_{\mathscr{D}^{\leq 0}(\mathcal{A})}\lra H^0:\;\;\mathscr{D}^{\leq 0}(\mathcal{A})\to \mathscr{D}^{\leq 0}(\mathcal{A}).$$
When $\mathcal{A}=A^{\rm e}\Modcat$ and $X, Y\in\mathscr{D}^{\leq 0}(\mathcal{A})$, we denote the composite of the following morphisms by
$$
\theta_{X,Y}:\xymatrix{X\otimesL_AY\ar[rr]^-{\xi_X\otimesL_A\xi_Y}&& H^0(X)\otimesL_AH^0(Y)\ar[rr]^-{\xi_{H^0(X)\otimesL_AH^0(Y)}}&&H^0(X)\otimes_AH^0(Y).}
$$
Then $\theta_{X,Y}$ is natural in $X$ and $Y$. This gives rise to a natural transformation
$$
\theta_{-,-}:\;\; (-)\otimesL_A(-)\lra H^0(-)\otimes_AH^0(-): \quad \mathscr{D}^{\leq 0}(\mathcal{A})\times \mathscr{D}^{\leq 0}(\mathcal{A})\lra A^{\rm e}\Modcat.
$$
We have the following result.

\begin{Lem} \label{IDH}
$(1)$ For $X\in \mathscr{D}^{\leq 0}(\mathcal{A})$, the morphism $H^0(\xi_X)$ is an automorphism of $H^0(X)$.

$(2)$ For a morphism $f:X\to Y$ in $\mathscr{D}^{\leq 0}(\mathcal{A})$, there is a unique morphism $f': H^0(X)\to H^0(Y)$ in $\mathcal{A}$ such that $f\xi_Y=\xi_X f'$. Moreover, $f'=H^0(\xi_X)^{-1}H^0(f)H^0(\xi_Y)$.

$(3)$ Let $\mathcal{A}:=A^{\rm e}\Modcat$. Then the map $H^0(\theta_{X,Y}): H^0(X\otimesL_AY)\to H^0(X)\otimes_AH^0(Y)$ is an isomorphism and $\theta_{X,Y}=\xi_{X\otimesL_AY}H^0(\theta_{X,Y})$. Thus there is a natural isomorphism of functors
$$H^0(\theta_{-,-}):\;\;H^0(-\otimesL_A-)\lraf{\simeq} H^0(-)\otimes_AH^0(-):\quad  \mathscr{D}^{\leq 0}(\mathcal{A})\times \mathscr{D}^{\leq 0}(\mathcal{A})\lra A^{\rm e}\Modcat.$$
\end{Lem}

{\it Proof.} $(1)$ and $(2)$ follow from the construction of $\xi$. Note that $H^0(\xi_X\otimesL_A\xi_Y)$ and $H^0(\xi_{H^0(X)\otimesL_AH^0(Y)})$ are isomorphisms.
Since $\xi_{H^0(X)\otimes_AH^0(Y)}$ is the identity, $(3)$ follows from $(2)$. $\square$

\medskip
In the rest of this section, let $\varphi:A\to A'$ be a homomorphism of algebras.
Define
$$W:=\Phi(A'), \;\;B':=H^0(W),\;\; W':=\tau^{\leq 0}W \; \; \mbox{and} \; \; \varphi':=H^0(\Phi(\varphi)): B\lra B'.$$

\begin{Lem}\label{alg-hom}
$(1)$ $A'\otimesL_A T$ is isomorphic in $\D{A'}$ to a tilting complex if and only if $H^n(W)=0$ for all $n\neq 0$.

$(2)$ $B'$ is an algebra with the multiplication induced from
$$\xymatrix{B'\otimes_BB'\;\;\ar[rr]^-{H^0(\theta_{W', W'})^{-1}}_-{\simeq}&&\;H^0(W'\otimesL_BW') \ar[rr]^-{H^0(\lambda_W^0\otimesL_B\lambda_W^0)}&& H^0(W\otimesL_BW)\ar[r]^-{H^0(\phi_{A',A'})}_-{\simeq}& H^0\big(\Phi(A'\otimesL_AA')\big)\ar[r]^-{H^0(\Phi(\pi))}& B'}$$
where  $\pi: A'\otimesL_AA'\to A'$ is the composite of $\xi_{A'\otimesL_AA'}: A'\otimesL_AA'\to A'\otimes_AA'$ with the multiplication map $A'\otimes_AA'\to A'$.

$(3)$ $B'$ and $\End_{\D{A'}}(A'\otimesL_A T)$ are isomorphic as algebras. Moreover, $\varphi'$ is a homomorphism of algebras.
\end{Lem}

{\it Proof.} $(1)$ Since $T$ is isomorphic in $\D{A}$ to a tilting complex $P$, we have $A'\otimesL_A T\simeq A'\otimes_AP$ in $\D{A'}$ and $A'\otimes_AP\in\Kb{\prj{A'}}$. As $\add(P)$ generates $\Kb{\prj{A}}$ as a triangulated category, $\add(A'\otimes_AP)$ generates $\Kb{\prj{A'}}$ as a triangulated category. This implies that $A'\otimesL_A T$ is isomorphic in $\D{A'}$ to a tilting complex if and only if $A'\otimesL_A T$ is self-orthogonal in $\D{A'}$.
Moreover, for $n\in\mathbb{Z}$, it follows from the isomorphism $\varepsilon_B:T^{\vee}\otimesL_AT\to B$ in $\D{B^{\rm e}}$ that there is a series of isomorphisms
$$(\ast)\;\;\Hom_{\D{A'}}(A'\otimesL_A T, A'\otimesL_A T[n])\simeq\Hom_{\D{A}}(T, A'\otimesL_A T[n])\simeq
\Hom_{\D{B}}(T^{\vee}\otimesL_AT, T^{\vee}\otimesL_AA'\otimesL_A T[n])$$
$$\quad\quad\quad\quad\quad\quad\quad\quad\quad\quad\quad\quad\quad
\simeq\Hom_{\D{B}}(B, T^{\vee}\otimesL_AA'\otimesL_A T[n])=\Hom_{\D{B}}(B, W[n])\simeq H^n(W).$$
Thus $A'\otimesL_A T$ is self-orthogonal if and only if $H^n(W)=0$ for all $n\neq 0$.
This shows $(1)$.

$(2)$ If taking $n=0$ in ($*$), we get an isomorphism $\End_{\D{A'}}(A'\otimesL_A T)\simeq B'$ of $k$-modules. Via the isomorphism, we can transfer the algebra structure of $\End_{\D{A'}}(A'\otimesL_A T)$ to the one of $B'$.

Let $s_i\in\Hom_{\D{B}}(B, W)$ for $i=1,2$. By $(\ast)$, there are morphisms $t_i: T\to  A'\otimesL_A T$ in $\D{A}$ such that $s_i=\varepsilon_B^{-1}(T^{\vee}\otimesL_At_i)$. By the first isomorphism in $(\ast)$, we  can define a multiplication on the abelian group $\Hom_{\D{A}}(T, A'\otimesL_A T)$, that is,
the multiplication of $t_1$ with $t_2$ is given by the composition of the morphisms
$$
t_1\cdot t_2:\;\; T\lraf{t_1}A'\otimesL_A T\lraf{A'\otimesL_At_2}A'\otimesL_A A'\otimesL_AT\lraf{\pi\otimesL_AT}A'\otimesL_A T.
$$
This yields the product $s_1\cdot s_2\in\Hom_{\D{B}}(B, W)$ of $s_1$ with $s_2$, described by the composite of the morphisms
$$
s_1\cdot s_2:\;\;B\lraf{\varepsilon_B^{-1}}T^{\vee}\otimesL_AT\lraf{T^{\vee}\otimesL_At_1\;}T^{\vee}\otimesL_AA'\otimesL_A T\lraf{T^{\vee}\otimesL_AA'\otimesL_At_2}T^{\vee}\otimesL_AA'\otimesL_A A'\otimesL_AT\lraf{\Phi(\pi)}T^{\vee}\otimesL_AA'\otimesL_A T=W.
$$
Since $({_A}T\otimesL_B-, {_B}T^{\vee}\otimesL_A-)$ is an adjoint pair of functors between $\D{B}$ and $\D{A}$, the composite
of the morphisms
$$
T\lraf{\rm can} T\otimesL_BB\lraf{T\otimesL_B\varepsilon_B^{-1}} T\otimesL_BT^{\vee}\otimesL_AT\lraf{\varepsilon_A\otimesL_AT}A\otimesL_AT\lraf{\rm{can}}T
$$
is the identity morphism of $T$, where the first and last morphisms are canonical isomorphisms. It follows that
$t_2$ is the composite of the morphisms
$$
{_A}T\lraf{\rm can} T\otimesL_BB\lraf{T\otimesL_B\varepsilon_B^{-1}} T\otimesL_BT^{\vee}\otimesL_AT\lraf{T\otimesL_BT^{\vee}\otimesL_At_2}T\otimesL_BT^{\vee}\otimesL_AA'\otimesL_AT\lraf{\varepsilon_A\otimesL_AA'\otimesL_AT}
A\otimesL_AA'\otimesL_AT\lraf{\rm{can}}{_A}A'\otimesL_AT,
$$
and therefore the multiplication $s_1\cdot s_2$ is the composite of the morphisms
$$
\xymatrix{B\ar[r]^-{s_1}& W\ar[r]^-{\rm can}_-{\simeq}
&W\otimesL_BB\ar[r]^-{W\otimesL_Bs_2}& W\otimesL_BW=\Phi(A')\otimesL_B\Phi(A')\ar[r]^-{\phi_{A',A'}}_-{\simeq}&\Phi(A'\otimesL_AA')\ar[r]^-{\Phi(\pi)}&W}.
$$
Since the inclusion $\lambda_W^0:W'\to W$ induces an isomorphism $\Hom_{\D{B}}(B,W')\simeq \Hom_{\D{B}}(B,W)$,
there are $s_i'\in\Hom_{\D{B}}(B,W')$ for $i=1,2$ such that $s_i=s_i'\lambda_W^0$.
Let $$\widetilde{s_i}:=s_i'\pi_{W'}^0\in \Hom_B(B,B').$$
Since $H^0$ induces an isomorphism $\Hom_{\D{B}}(B,W)\simeq \Hom_B(B, B')\simeq B'$ as $k$-modules, we have $H^0(s_i)=\widetilde{s_i}$.
In the diagram
$$
\xymatrix{B\ar[r]^-{s_1}\ar@{=}[d]& W\ar[r]^-{\rm can}_-{\simeq}&W\otimesL_BB\ar[r]^-{W\otimesL_Bs_2}& W\otimesL_BW\ar[r]^-{{\phi_{A',A'}}\;\Phi(\pi)}& W\\
B\ar[r]^-{s_1'}\ar@{=}[d]&W'\ar[r]^-{\rm can}_-{\simeq}\ar[u]^-{\lambda_W^0}\ar[d]_-{\pi_{W'}^0}& W'\otimesL_BB
\ar[r]^-{W'\otimesL_Bs_2'}\ar[u]^-{\lambda_W^0\otimesL_BB}\ar[d]_-{\pi_{W'}^0\otimesL_BB}&\;W'\otimesL_BW'\ar[u]_-{\lambda_W^0\otimesL_B\lambda_W^0}
\ar[d]^-{\pi_{W'}^0\otimesL_B\pi_{W'}^0}\ar@{=}[r]&W'\otimesL_BW'\ar[dd]^-{\theta_{W',W'}}\ar[u]\\
B\ar@{=}[d]\ar[r]^-{\widetilde{s_1}}&B'\ar@{=}[d]\ar[r]^-{\rm can}_-{\simeq}& B'\otimesL_BB\ar[d]_-{\pi^0_{B'\otimesL_BB}}
\ar[r]^-{B'\otimesL_B\widetilde{s_2}}&\;B'\otimesL_BB'\ar[d]^-{\pi^0_{B'\otimesL_BB'}}&\\
B\ar[r]^-{\widetilde{s_1}}&B'\ar[r]^-{\rm can}_-{\simeq}& B'\otimes_BB
\ar[r]^-{B'\otimes_B\widetilde{s_2}}&\;B'\otimes_BB'\ar@{=}[r]&B'\otimes_BB'}
$$
of morphisms in $\D{B}$, all the squares are commutative and $H^0(\theta_{W',W'})$ is an isomorphism. Let $\mu: B\to B'\otimes_BB'$ be the composite of the morphisms in the bottom line of the diagram. Then
$$H^0(s_1\cdot s_2)=\mu\, H^0(\theta_{W',W'})^{-1}H^0(\lambda_W^0\otimesL_B\lambda_W^0)H^0(\phi_{A',A'})H^0(\Phi(\pi)):\;B\lra B'.$$
Note that $\mu$ sends the identity $1$ of $B$ to $(1)\widetilde{s_1} \otimes (1)\widetilde{s_2}$. Now, by identifying $B'$ with $\Hom_{\D{B}}(B,W)$ and also with $\Hom_B(B,B')$,  $(2)$ can be proved.

$(3)$ The isomorphism $\End_{\D{A'}}(A'\otimesL_A T)\simeq B'$ as algebras has been shown in the proof of $(2)$. Now, we denote by $\mu_{B'}:B'\otimes_BB'\to B'$ the composite of the morphisms in $(2)$. Recall that $\Phi(A)\simeq B$ and $H^0(\Phi(A))=B$. If $A'=A$ and $\varphi={\rm Id}_A$, then $B'=B$ and $\mu_B: B\otimes_BB\to B$ is the canonical isomorphism induced by the multiplication of $B$. For a general $\varphi:A\to A'$, there is an equality $\mu_B \varphi'=(\varphi'\otimes_B\varphi')\mu_{B'}$ which means that $\varphi'$ is an algebra homomorphism. Note that if $B$ and $B'$ are identified with $\End_{\D{A}}(T)$ and $\End_{\D{A'}}(A'\otimesL_A T)$, respectively, then $\varphi':B\to B'$ is exactly the algebra homomorphism induced from $A'\otimesL_A-:\D{A}\to\D{A'}$. $\square$

\medskip
The following result provides a method for constructing derived equivalences of algebras with idempotents from
given ones. It also generalizes derived equivalences of
trivial extensions of algebras by bimodules (see \cite[Corollary 5.4]{Rickard3}).

\begin{Prop}\label{RDE}
Suppose $H^n(W)=0$ for all $n\neq 0$. Then the pairs $(A', (e)\varphi)$ and $(B', (f)\varphi')$ of algebras with idempotents are derived equivalent. In particular, $A'$ and $B'$ are derived equivalent.
\end{Prop}

{\it Proof.} By Lemma \ref{alg-hom}, $A'\otimesL_A T$ is isomorphic in $\D{A'}$ to a tilting complex
and $\End_{\D{A'}}(A'\otimesL_A T)\simeq B'$ as algebras. It follows from Theorem \ref{Derived equivalence} that $A'$ and $B'$ are derived
equivalent. Clearly,  $T\otimesL_BBf\simeq Tf\in\add({_A}T)$ and $A'\otimesL_ATf\in\add({_{A'}}A'\otimesL_A T)$.
Since $(A, e)$ and $(B,f)$ are derived equivalent, $Tf$ generates $\Kb{\add(Ae)}$ in $\D{A}$ by Lemma \ref{IMDE}(2).
Applying the functor $A'\otimesL_A-:\D{A}\to\D{A'}$ to $Tf$ and $\Kb{\add(Ae)}$, we see from  the general result in the proof of Lemma \ref{EVA}(2)
that $A'\otimesL_ATf$ generates $\Kb{\add(A'(e)\varphi)}$ in $\D{A'}$.
Note that $$\Hom_{\D{A'}}(A'\otimesL_AT, A'\otimesL_ATf)\simeq\Hom_{\D{A}}(T, A'\otimesL_A Tf)\simeq
\Hom_{\D{B}}(T^{\vee}\otimesL_AT, T^{\vee}\otimesL_AA'\otimesL_A Tf)$$
$$\quad\quad\quad\quad\quad\quad\quad
\simeq\Hom_{\D{B}}(B, W\otimes_BBf)\simeq H^0(W)\otimes_BBf\simeq B'(f)\varphi'.$$
By Lemma \ref{IMDE}(2), the pairs $(A', (e)\varphi)$ and $(B', (f)\varphi')$ are derived equivalent. $\square$

\medskip
Now, we turn to mirror-reflective algebras at any levels. Recall that, for each $\lambda\in Z(\Lambda)$, the  multiplication map $(\cdot\lambda): \Lambda\to\Lambda$ induces a homomorphism
$\omega_\lambda: \Delta_0\otimes_A \Delta_0\to \Delta_0$ in $A^{\rm e}\Modcat$, which is the composite of the maps
$$
\Delta_0\otimes_A\Delta_0\lraf{\omega_e}\Delta_0\lraf{F_e(\cdot\lambda)}
\Delta_0.
$$
We define the derived version of $\omega_\lambda$ to be the composite of the maps in $\D{A^{\rm e}}$:
$$
\mathbb{L}{\omega_\lambda}:\;\;
\Delta\otimesL_A\Delta\lraf{\simeq}\Delta\lraf{\mathbb{L}F_e(\cdot\lambda)}
\Delta
$$
where the first isomorphism is canonical, due to $eA\otimesL_AAe\simeq\Lambda$ in $\D{\Lambda^{\rm e}}$. Note that
both $\omega_e$ and $\mathbb{L}{\omega_e}$ are isomorphisms since $(\cdot e)$ is the identity map of $\Lambda$, and that $F_e$ and $\mathbb{L}F_e$ are fully faithful functors. Thus we have the following result.

\begin{Lem}\label{Isomorphism}
There are isomorphisms
$$\omega_{(-)}:\;\;Z(\Lambda)\lraf{\simeq} \Hom_{A^{\rm e}}(\Delta_0\otimes_A \Delta_0, \Delta_0),\;\;\lambda\mapsto \omega_\lambda=\omega_eF_e(\cdot\lambda),$$
$$\mathbb{L}\omega_{(-)}:\;\;Z(\Lambda)\lraf{\simeq} \Hom_{\D{A^{\rm e}}}(\Delta\otimesL_A\Delta,\Delta),\;\; \lambda\mapsto \mathbb{L}\omega_\lambda=(\mathbb{L}\omega_e)\mathbb{L}F_e(\cdot\lambda).$$
Moreover, $\omega_\lambda$ is an isomorphism if and only if $\lambda$ is invertible if and only if $\mathbb{L}\omega_\lambda$ is an isomorphism.
\end{Lem}

Similarly, for $\mu\in Z(\Gamma)$, there is a homomorphism $\omega_{\mu}: \Theta_0\otimes_A \Theta_0\to \Theta_0$ in $B^{\rm e}\Modcat$ with its derived version
$$\mathbb{L}\omega_{\mu}:\;\;
\Theta\otimesL_B\Theta\lraf{\simeq}\Theta
\lraf{\mathbb{L}F_f(\cdot\mu)}\Theta
$$
in $\D{B^{\rm e}}$.
Following the diagram $(\sharp)$, let $\eta:\Phi\circ \mathbb{L}F_e\to \mathbb{L}F_f\circ \Phi'$ be a natural isomorphism of functors from $\D{\Lambda^{\rm e}}$ to $\D{B^{\rm e}}$ and let $\tau:\Phi'(\Lambda)\to \Gamma$ be an isomorphism in $\D{\Gamma^{\rm e}}$.

\begin{Lem}\label{CM}
The following hold for $\lambda\in Z(\Lambda)$ and $\mu:=(\lambda)\sigma\in Z(\Gamma)$.

$(1)$ There are commutative diagrams
$$\xymatrix{\Delta\otimesL_A\Delta\ar[r]^-{\mathbb{L}{\omega_\lambda}}
\ar[d]_-{\theta_{\Delta,\Delta}}&\Delta\ar[d]^-{\xi_{\Delta}}\\
\Delta_0\otimes_A\Delta_0\ar[r]^-{\omega_{\lambda\lambda_1}}&\Delta_0 ,}
\quad\quad \xymatrix{\Phi(\Delta)\otimesL_B\Phi(\Delta)\ar[rr]^-{\phi_{\Delta,\Delta}\,\Phi(\mathbb{L}{\omega_\lambda})}
\ar[d]_-{\tau_1\otimesL_B\tau_1}&&\Phi(\Delta)\ar[d]^-{\tau_1}\\
\Theta\otimesL_B\Theta\ar[rr]^-{\mathbb{L}\omega_{\mu\mu_1}}&& \Theta}$$
where $\lambda_1\in Z(\Lambda)$ and $\mu_1\in Z(\Gamma)$ are invertible, and $\tau_1:=\eta_\Lambda\mathbb{L}F_f(\tau):\Phi(\Delta)\to \Theta$ is an isomorphism.

$(2)$ Let $$W_0:=\Phi(\Delta_0),\;\;\tau_2:=\tau_1^{-1}\Phi(\xi_\Delta):\;\;\Theta\lra W_0,\;\; \psi_\lambda:=\xi_{\Delta_0\otimesL_A\Delta_0}\omega_{\lambda\lambda_1}:\;\;\Delta_0\otimesL_A\Delta_0\lra \Delta_0.$$
Then there is a commutative diagram in $\D{B^{\rm {e}}}:$
$$\xymatrix{\Theta\otimesL_B\Theta\ar[d]_-{\tau_2\otimesL_B\tau_2}\ar[rrr]^-{\mathbb{L}\omega_{\mu\mu_1}}&&& \Theta\ar[d]^-{\tau_2}\\
W_0\otimesL_BW_0\ar[rrr]^-{\phi_{\Delta_0,\Delta_0}\,\Phi(\psi_\lambda)}
&&&W_0.}$$

$(3)$ If $W_0$ lies in $\mathscr{D}^{\leq 0}(B^{\rm e})$, then there is a commutative diagram in $B^{\rm e}\Modcat:$
$$\xymatrix{\Theta_0\otimes_B\Theta_0\ar[d]_-{H^0(\tau_2)\otimes_B H^0(\tau_2)}\ar[rrrr]^-{\omega_{\mu\mu\,'}}&&&& \Theta_0\ar[d]^-{H^0(\tau_2)}\\
H^0(W_0)\otimes_BH^0(W_0)\;\ar[rrrr]^-{H^0(\theta_{\,W_0,W_0})^{-1}
H^0(\phi_{\Delta_0,\Delta_0})\,H^0(\Phi(\psi_\lambda))}&&&&\;H^0(W_0),}$$
where $\mu\,'\in Z(\Gamma)$ is invertible. If, in addition, $H^0(\Phi(\xi_\Delta))$ is an isomorphism, then there is an algebra isomorphism
$$H^0\big(\Phi(R(A,e,\lambda))\big)\simeq R(B,f,\mu).$$

$(4)$ If $W_0\simeq \Theta_0$ in $\D{B^{\rm e}}$, then $H^0(\Phi(\xi_\Delta))$ is an isomorphism of $B^e$-modules.
\end{Lem}

{\it Proof.} $(1)$ Note that $\mathbb{L}{\omega_e}$ is an isomorphism. By Lemma \ref{IDH}(2)-(3), $\mathbb{L}F_e(\cdot \lambda)\xi_\Delta=\xi_\Delta F_e(\cdot\lambda)$ and there is a unique isomorphism $\alpha:\Delta_0\otimes_A\Delta_0\to \Delta_0$ such that
$(\mathbb{L}{\omega_e})\xi_{\Delta}=\theta_{\Delta,\Delta}\alpha$. By the first isomorphism in Lemma \ref{Isomorphism}, there is an element $\lambda_1\in Z(\Lambda)^{\times}$ such that $\alpha=\omega_{\lambda_1}=\omega_eF_e(\cdot\lambda_1).$ Thus
$$(\mathbb{L}{\omega_\lambda})\xi_\Delta=(\mathbb{L}{\omega_e})\mathbb{L}F_e(\cdot \lambda)\xi_\Delta=\theta_{\Delta,\Delta}\alpha F_e(\cdot\lambda)=\theta_{\Delta,\Delta}\omega_eF_e(\cdot\lambda_1)F_e(\cdot\lambda)=\theta_{\Delta,\Delta}\omega_eF_e(\cdot(\lambda\lambda_1))
=\theta_{\Delta,\Delta}\omega_{\lambda\lambda_1}.$$ Hence the diagram in the left-hand side of (1) is commutative.

Since $\mathbb{L}{\omega_e}$ and $\tau$ are isomorphisms and $\eta$ is a natural isomorphism, there is a unique isomorphism $\beta:\Theta\otimesL_B\Theta\to\Theta$ such that all the squares in the diagram are commutative:
$$
\xymatrix{\Phi(\Delta)\otimesL_B\Phi(\Delta)\ar[r]^-{\phi_{\Delta, \Delta}}\ar[d]_-{\eta_\Lambda\otimesL_B\eta_\Lambda}&\Phi(\Delta\otimesL_A\Delta)
\ar[r]^-{\Phi(\mathbb{L}{\omega_e})}&\Phi(\Delta)\ar[d]_-{\eta_\Lambda}\ar[rr]^-{\Phi\circ \mathbb{L}F_e(\cdot\lambda)}&&\Phi(\Delta)\ar[d]^-{\eta_\Lambda}\\
\mathbb{L}F_f\circ\Phi'(\Lambda)\otimesL_B\mathbb{L}F_f\circ \Phi'(\Lambda)\ar[d]_-{\mathbb{L}F_f(\tau)\otimesL_B\mathbb{L}F_f(\tau)}&&\mathbb{L}F_f\circ \Phi'(\Lambda)\ar[d]_-{\mathbb{L}F_f(\tau)}\ar[rr]^-{\mathbb{L}F_f\circ \Phi'(\cdot\lambda)}&&\mathbb{L}F_f\circ\Phi'(\Lambda)\ar[d]^-{\mathbb{L}F_f(\tau)}\\
\Theta\otimesL_B\Theta\ar[rr]^-{\beta}&&\Theta\ar[rr]^-{\mathbb{L}F_f(\cdot\mu)}&&\Theta}
$$
Further, by applying a similar isomorphism in Lemma \ref{Isomorphism} to the pair $(B,f)$, we have $\beta=\mathbb{L}\omega_{\mu_1}=(\mathbb{L}\omega_f)\mathbb{L}F_f(\cdot\mu_1)$ for some invertible element $\mu_1\in Z(\Gamma)$.
This implies $$\beta \mathbb{L}F_f(\cdot\mu)=(\mathbb{L}\omega_f)\mathbb{L}F_f(\cdot\mu_1)\mathbb{L}F_f(\cdot\mu)
=(\mathbb{L}\omega_f)\mathbb{L}F_f(\cdot(\mu\mu_1))=\mathbb{L}\omega_{\mu\mu_1}.$$
Thus the second diagram in (1) is commutative.

$(2)$ It follows from $\theta_{\Delta,\Delta}=(\xi_\Delta\otimesL_A\xi_\Delta)\xi_{\Delta_0\otimesL_A\Delta_0}$ and the first diagram in (1) that
there is the commutative diagram:
$$\xymatrix{\Delta\otimesL_A\Delta\ar[r]^-{\mathbb{L}{\omega_\lambda}}
\ar[d]_-{\xi_\Delta\otimesL_A\xi_\Delta}&\Delta\ar[d]^-{\xi_{\Delta}}\\
\Delta_0\otimesL_A\Delta_0\ar[r]^-{\psi_{\lambda}}&\Delta_0.}$$
Applying $\Phi$ and the natural isomorphism $\phi_{-,-}:\;\; \Phi(-)\otimesL_B\Phi(-)\lraf{\simeq}\Phi(-\otimesL_A-)$ to this diagram,
we get another commutative diagram:
$$\xymatrix{\Phi(\Delta)\otimesL_B\Phi(\Delta)\ar[rr]^-{\phi_{\Delta,\Delta}\,\Phi(\mathbb{L}{\omega_\lambda})}
\ar[d]_-{\Phi(\xi_\Delta)\otimesL_B\Phi(\xi_\Delta)}&&\Phi(\Delta)\ar[d]^-{\Phi(\xi_{\Delta})}\\
W_0\otimesL_AW_0\ar[rr]^-{\phi_{\Delta_0, \Delta_0}\Phi(\psi_{\lambda})}&&W_0.}$$
Now, the commutative diagram in $(2)$ follows from the second commutative diagram in $(1)$.

$(3)$ Note that $H^0\circ \mathbb{L}F_e=F_e$. Applying $H^0$ to the diagram in $(2)$, we see from Lemma \ref{IDH}(3) that the squares in the  diagram
$$
(\natural_1)\quad \xymatrix{\Theta_0\otimes_B\Theta_0\ar[d]_-{H^0(\tau_2)\otimes_BH^0(\tau_2)}&&\ar[ll]_-{H^0(\theta_{\Theta, \Theta})}^-{\simeq}
H^0(\Theta\otimesL_B\Theta)\ar[d]_-{H^0(\tau_2\otimesL_B\tau_2)}\ar[rrr]^-{H^0(\mathbb{L}\omega_{\mu\mu_1})}&&& H^0(\Theta)\ar[d]^-{H^0(\tau_2)}\\
H^0(W_0)\otimes_BH^0(W_0)&&\ar[ll]_-{H^0(\theta_{W_0, W_0})}^-{\simeq} H^0(W_0\otimesL_BW_0)\ar[rrr]^-{H^0(\phi_{\Delta_0,\Delta_0})\,H^0(\Phi(\psi_\lambda))}&&& H^0(W_0)}
$$
are commutative, where the isomorphisms are due to $\Theta, W_0\in\mathscr{D}^{\leq 0}(B^{\rm e})$.
Moreover, for the pair $(B,f)$ and $\mu\mu_1\in Z(\Gamma)$,  we obtain similarly the following commutative
diagrams, in which the second one is obtained from the first by the functor $H^0$:
$$
(\natural_2)\quad
\xymatrix{\Theta\otimesL_A\Theta\ar[r]^-{\mathbb{L}{\omega_{\mu\mu_1}}}
\ar[d]_-{\theta_{\Theta,\Theta}}&\Theta\ar[d]^-{\xi_{\Theta}}\\
\Theta_0\otimes_A\Theta_0\ar[r]^-{\omega_{\mu\mu_1\mu_2}}&\Theta_0 \;,}\qquad
\xymatrix{H^0(\Theta\otimesL_A\Theta)\ar[rr]^-{H^0(\mathbb{L}{\omega_{\mu\mu_1}})}
\ar[d]_-{H^0(\theta_{\Theta,\Theta})}&&\Theta_0\ar[d]^-{H^0(\xi_{\Theta})}\\
\Theta_0\otimes_A\Theta_0\ar[rr]^-{\omega_{\mu\mu_1\mu_2}}&&\Theta_0}$$
where $\mu_2\in Z(\Gamma)$ is invertible and $H^0(\xi_{\Theta})$ is an automorphism by Lemma \ref{IDH}(1). Since the functor $F_f$ induces an algebra isomorphism $Z(\Gamma)\to\End_{B^{\rm e}}(\Theta_0)$, there is an invertible element $\mu_3\in Z(\Gamma)$ such that $H^0(\xi_{\Theta})=F_f(\cdot \mu_3)$. This implies
$$
(\natural_3)\quad\omega_{\mu\mu_1\mu_2}H^0(\xi_{\Theta})^{-1}=\omega_{\mu\mu_1\mu_2}F_f(\cdot \mu_3^{-1})=\omega_{\mu\mu_1\mu_2\mu_3^{-1}}.
$$
Let $\mu\,':=\mu_1\mu_2\mu_3^{-1}\in Z(\Gamma)$. Then $\mu'$ is invertible. By $(\natural_1)$-$(\natural_3)$, we obtain the commutative diagram in $(3)$.

Next, we apply Lemma \ref{alg-hom}(2) to show the algebra isomorphism in $(3)$.

Let $A':=R(A,e,\lambda\lambda_1)$, $B':=H^0(\Phi(A'))$ and $\varphi:A\to A'$ be the canonical injection. By Lemma \ref{alg-hom}(2), $B'$ is an algebra.
Since $A'=A\oplus\Delta_0$ and $\Phi(A)\simeq B$,  there holds $\Phi(A')\simeq B\oplus W_0$ in $\D{B^{\rm e}}$. Now, we identity $B'$ with $B\oplus H^0(W_0)$ as $B^{\rm e}$-modules and describe the multiplication of $B'$ in terms of the one of $A'$ and the one in Lemma \ref{alg-hom}(2):

The multiplication of $B$ with $B'$ (or $B'$ with $B$) is given by left (or right) multiplication
since $B'$ is a $B$-$B$-bimodule; while the multiplication on $H^0(W_0)$  is induced from the composition
$$\xymatrix{H^0(W_0)\otimes_BH^0(W_0) \ar[rr]^-{H^0(\theta_{W_0', W_0'})^{-1}}&&H^0(W_0'\otimesL_BW_0') \ar[rr]^-{H^0(\lambda_{W_0}^0\otimesL_B\lambda_{W_0}^0)}&& H^0(W_0\otimesL_BW_0)\ar[d]^-{H^0(\phi_{\Delta_0, \Delta_0})}\\
H^0(W_0)&&H^0\big(\Phi(\Delta_0\otimes_A\Delta_0)\big)\ar[ll]_-{H^0(\Phi(\omega_{\lambda\lambda_1}))}&& H^0\big(\Phi(\Delta_0\otimesL_A\Delta_0)\big)\ar[ll]_-{H^0(\Phi(\xi_{\Delta_0\otimesL_A\Delta_0}))}}$$
where $W_0':=\tau^{\leq 0}W_0$ and the injection $\lambda_{W_0}^0: W_0'\to W_0$ is an isomorphism by $W_0\in\mathscr{D}^{\leq 0}(B^{\rm e})$.
It then follows from $$\theta_{W_0, W_0}=\big((\lambda_{W_0}^0)^{-1}\otimesL_B(\lambda_{W_0}^0)^{-1}\big)\theta_{W_0', W_0'}=
\big(\lambda_{W_0}^0\otimesL_B\lambda_{W_0}^0\big)^{-1}\theta_{W_0', W_0'},$$
that $H^0(\theta_{\,W_0,W_0})^{-1}=H^0(\theta_{W_0', W_0'})^{-1}H^0(\lambda_{W_0}^0\otimesL_B\lambda_{W_0}^0).$
Thus the multiplication of $H^0(W_0)$ with $H^0(W_0)$ in $B'$ is induced from
$$H^0(\theta_{\,W_0,W_0})^{-1}H^0(\phi_{\Delta_0,\Delta_0})\,H^0(\Phi(\psi_\lambda)):\;\; H^0(W_0)\otimes_BH^0(W_0) \lra H^0(W_0).$$

Suppose that $H^0(\Phi(\xi_\Delta))$ is an isomorphism. Then $H^0(\tau_2)$ is an isomorphism and $B'\simeq B\oplus \Theta_0$ as $B^{\rm e}$-modules. Moreover, the commutative diagram in $(3)$ implies that $H^0(\tau_2)$ induces an algebra isomorphism $R(B, f, \mu\mu\,')\simeq B'$ which lifts the identity map of $B$. Since $\lambda_1\in Z(\Lambda)$ and $\mu\,'\in Z(\Gamma)$ are invertible, it follows from \cite[Lemma 3.2(2)]{xcf} that
$A'\simeq R(A,e,\lambda)$ and $R(B, f, \mu\mu\,')\simeq R(B, f, \mu)$ as algebras.
Thus there are algebra isomorphisms $H^0\big(\Phi(R(A,e,\lambda))\big)\simeq H^0\big(\Phi(A')\big)=B'\simeq R(B,f,\mu)$.

$(4)$ Under the identifications $G_e(\Delta_0)=\Lambda$ and $\Delta=\mathbb{L}F_e(\Lambda)$, we see that
$\xi_\Delta:\Delta=(\mathbb{L}F_e\circ G_e)(\Delta_0)\to\Delta_0$ is the counit adjunction morphism of $\Delta_0$
associated with the adjoint pair $(\mathbb{L}F_e, G_e)$. Similarly, up to isomorphism, $\xi_\Theta:\Theta=(\mathbb{L}F_f\circ G_f)(\Theta_0)\to\Theta_0$
is the counit adjunction morphism of $\Theta_0$ associated with the adjoint pair $(\mathbb{L}F_f, G_f)$.
Now, recall that two morphisms $f_i: X_i\to Y_i$ for $i=1,2$ in an additive category
are isomorphic if there are isomorphisms $\alpha_1: X_1\to X_2$ and $\alpha_2:Y_1\to Y_2$ such that $f_1\alpha_2=\alpha_1 f_2$.
By the diagram $(\sharp)$, the functor $\Phi$ is an equivalence and there is a natural isomorphism
$$\Phi\circ \mathbb{L}F_e\circ G_e\lraf{\simeq} \mathbb{L}F_f\circ G_f\circ \Phi:\;\;\D{A^{\rm e}}\lra \D{B^{\rm e}}.$$
This implies that $\Phi(\xi_\Delta):\Phi(\Delta)\to W_0$ is isomorphic to the counit adjunction morphism of $W_0$
associated with $(\mathbb{L}F_f, G_f)$. If $W_0\simeq \Theta_0$ in $\D{B^{\rm e}}$, then $\xi_\Theta$ and $\Phi(\xi_\Delta)$ are isomorphic as morphisms in $\D{B^{\rm e}}$. Since $H^0(\xi_\Theta)$ is an isomorphism by Lemma \ref{IDH}(1), $H^0(\Phi(\xi_\Delta))$ is an isomorphism. This shows $(4)$.
$\square$

\medskip
{\bf Proof of Theorem \ref{derived equivalent}}. For  each $\lambda\in Z(\Lambda)$, let $A':=R(A,e,\lambda)$, $\varphi:A\to A'$ the canonical injection and $B':=H^0(\Phi(A'))$. Since  $A'=A\oplus\Delta_0$ and $\Phi(A)\simeq B$, we have $\Phi(A')\simeq B\oplus \Phi(\Delta_0)$.
By assumption, $\Phi(\Delta_0)\simeq \Theta_0$ in $\D{B^{\rm e}}$. This implies $\Phi(A')\simeq B\oplus \Theta_0$ in $\D{B^{\rm e}}$, and therefore
$B'=B\oplus H^0(\Phi(\Delta_0))\simeq B\oplus \Theta_0$ and  $H^n(\Phi(A'))=0$ for all $n\neq 0$. Now, let $\varphi':=H^0(\Phi(\varphi)): B\to B'$.
By the multiplication of $B'$ in Lemma \ref{alg-hom}(2), $\varphi'$ is the canonical injection. Then $(e)\varphi=e\in A'$ and $(f)\varphi'=f\in B'$. By Proposition \ref{RDE}, $(A', e)$ and $(B', f)$ are derived equivalent.  Since $\Phi(\Delta_0)\simeq \Theta_0$ in $\D{B^{\rm e}}$, it follows from Lemma \ref{CM}(3)(4) that there is an algebra isomorphism $B'\simeq R(B,f, (\lambda)\sigma)$ which lifts the identity map of $B$. Consequently, $(A', e)$ and $\big(R(B,f,(\lambda)\sigma), f\big)$ are derived equivalent. Clearly, $(e)\sigma=f$ since $e$ and $f$ are identities of $\Lambda$ and $\Gamma$, respectively. Thus
$\big(R(A,e), e\big)$ and $\big(R(B,f), f\big)$ are derived equivalent. $\square$

\smallskip
A sufficient condition for the isomorphism in Theorem \ref{derived equivalent} to hold true is the vanishing of positive Tor-groups over corner algebras.

\begin{Prop}\label{VTG}
Suppose that there is a derived equivalence between $(A, e)$ and $(B,f)$ of algebras with idempotents, which is induced by a two-sided tilting complex ${_A}T_{B}$. If $\Tor_n^\Lambda(Ae, eA)=0=\Tor_n^\Gamma(Bf, fB)$ for all $n\geq 1$, then the derived equivalence $\Phi:\D{A^{\rm e}}\to \D{B^{\rm e}}$ associated with $T$ between the enveloping algebras $A^{\rm e}$ and $B^{\rm e}$ satisfies
$\Phi(Ae\otimes_\Lambda eA)\simeq Bf\otimes_\Gamma fB$ in $\D{B^{\rm e}}$.
\end{Prop}

{\it Proof.} Since $\Tor_n^\Lambda(Ae, eA)=0$ for all $n\geq 1$, we have $Ae\otimesL_\Lambda eA\simeq Ae\otimes_\Lambda eA$ in $\D{A^{\rm e}}$.
Similarly, $Bf\otimesL_\Gamma fB\simeq  Bf\otimes_\Gamma fB$ in $\D{B^{\rm e}}$. Moreover, since $(A, e)$ and $(B,f)$ are derived equivalent,
it follows from Lemma \ref{IMDE}(4) that $\Phi(Ae\otimesL_\Lambda eA)\simeq Bf\otimesL_\Gamma fB$ in $\D{B^{\rm e}}$ .
Thus $\Phi(Ae\otimes_\Lambda eA)\simeq Bf\otimes_\Gamma fB$.
$\square$

\medskip
{\bf Proof of Theorem \ref{main result}}.
Suppose that $A$ and $B$ are derived equivalent, gendo-symmetric algebras. Then the pair $(A,e)$ and $(B,f)$ are derived equivalent by Proposition \ref{DGS}. Without loss of generality, we assume that the derived equivalence between $(A, e)$ and $(B,f)$ is induced by a two-sided tilting complex $T\in\D{A\otimes_kB\opp}$. This gives rise to a derived equivalence between $A^{\rm e}$ and $B^{\rm e}$. Let
$\Phi:=T^\vee\otimesL_A-\otimesL_AT:\D{A^{\rm e}}\to \D{B^{\rm e}}$ be the associated equivalence.
Then $\Phi$ induces an algebra isomorphism $\sigma:Z(eAe)\to Z(fBf)$ (see the lines just before Theorem \ref{derived equivalent}).
Note that, for the gendo-symmetric algebra $(A,e)$, there is an isomorphism ${_A}Ae\otimes_\Lambda eA{_A}\simeq D(A)$ of $A$-$A$-bimodules
by \cite[Section 2.2]{FK16} or \cite[Lemma 4.1(2)]{xcf}. Similarly, ${_B}Bf\otimes_\Gamma fB{_B}\simeq D(B)$ as $B$-$B$-bimodules.
Since $\Phi(D(A))\simeq D(B)$  in $\D{B^{\rm e}}$, we have $\Phi(Ae\otimes_\Lambda eA)\simeq Bf\otimes_\Gamma fB$ in $\D{B^{\rm e}}$.
Now, Theorem \ref{main result} follows immediately from Theorem \ref{derived equivalent}. $\square$

\medskip
Finally, we present an example to illustrate the main result.
We consider the truncated polynomial algebra $\Lambda:=k[x]/(x^3)$. Let $X$ be the simple $\Lambda$-module and $Y$  the indecomposable $\Lambda$-module of length $2$. Then $A:=\End_{\Lambda}(\Lambda\oplus X)$ and $B:=\End_{\Lambda}(\Lambda\oplus Y)$ are derived equivalent, gendo-symmetric algebras. In this case, $Ae=\Hom_{\Lambda}(\Lambda\oplus X, \Lambda)$ and $Bf=\Hom_{\Lambda}(\Lambda\oplus Y,\Lambda)$. Clearly, $eAe\simeq eBe\simeq \Lambda.$ Moreover, $A$ and $B$ are given by the following quivers with relations, respectively:
$$
\xymatrix{
A:\quad &1\bullet\ar@(lu,ld)_{\gamma\;}\ar^-{\beta}@/^0.6pc/[r]
&\bullet 2,\ar^-{\alpha}@/^0.6pc/[l]&\quad
\qquad\qquad\quad B: & 1\bullet\ar^-{\beta}@/^0.6pc/[r]&\bullet 2,\ar^-{\alpha}@/^0.6pc/[l]}$$
$$ \qquad\qquad \alpha\beta=\alpha\gamma=0,\;\gamma^2=\beta\alpha.        \quad    \hspace{2.5cm} \quad\qquad\quad \alpha\beta\alpha\beta=0.\qquad $$
Further, $e$ and $f$ are corresponding to the vertex $1$ in the quivers, respectively. $R(A,e) $ and $R(B,f)$ are presented by the following quivers with relations, respectively.
$$
\xymatrix{
R(A, e):\quad &1\bullet\ar@(lu,ld)_{\gamma\;}\ar^-{\beta}@/^0.6pc/[r]
&\bullet 2\ar^-{\alpha}@/^0.6pc/[l]\ar_-{\bar{\alpha}}@/_0.6pc/[r]
&\bullet\bar{1}\;\ar@(ur,dr)^-{\bar{\gamma}\;,}\ar_-{\bar{\beta}}@/_0.6pc/[l]&\qquad
R(B,f):& 1\bullet\ar^-{\beta}@/^0.6pc/[r]&\bullet 2\ar^-{\alpha}@/^0.6pc/[l]\ar_-{\bar{\alpha}}@/_0.6pc/[r]&\bullet \bar{1},\ar_-{\bar{\beta}}@/_0.6pc/[l],}\qquad$$
$$ \qquad\qquad \beta\bar{\alpha}=\bar{\beta}\alpha=0,\,\;\alpha\gamma=\bar{\alpha}\bar{\gamma}=0,        \quad  \hspace{2.5cm} \quad\beta\bar{\alpha}=\bar{\beta}\alpha=0,\qquad$$
$$\qquad\quad
\gamma^2=\beta\alpha,\;\,\bar{\gamma}^2=\bar{\beta}\bar{\alpha},\;
\alpha\beta+\bar{\alpha}\bar{\beta}=0. \quad   \hspace{1.5cm} \quad\quad\quad\;\alpha\beta\alpha\beta+\bar{\alpha}\bar{\beta}\bar{\alpha}\bar{\beta}=0.
$$
By Theorem \ref{main result}, $R(A,e) $ and $R(B,f)$ are derived equivalent.

{\footnotesize
\smallskip
Hongxing Chen,

School of Mathematical Sciences  \&  Academy for Multidisciplinary Studies, Capital Normal University, 100048
Beijing, P. R. China;

{\tt Email: chenhx@cnu.edu.cn (H.X.Chen)}

\smallskip
Changchang Xi,

School of Mathematical Sciences, Capital Normal University, 100048
Beijing; \&
School of Mathematics and Statistics, Shaanxi Normal University, 710119 Xi'an, P. R. China

{\tt Email: xicc@cnu.edu.cn (C.C.Xi)}

\end{document}